# INFINITE DIMENSIONAL FAMILIES OF LOCALLY NONSOLVABLE PARTIAL DIFFERENTIAL OPERATORS

MICHAEL CHRIST, G. E. KARADZHOV, AND DETLEF MÜLLER

ABSTRACT. Local solvability is analyzed for natural families of partial differential operators having double characteristics. In some families the set of all operators that are not locally solvable is shown to have both infinite dimension and infinite codimension.

## 1. INTRODUCTION

A linear partial differential operator $L$ is said to be locally solvable at a point $p$ if there exists an open neighborhood $V$ of $p$ such that for any $f \in C_0^\infty(V)$ there exists $u \in \mathcal{D}'(\mathcal{V})$ satisfying $Lu = f$ in $V$. Whereas there exists a rather complete theory for operators of principal type with $C^\infty$ coefficients [1],[11],[7], the theory for operators with multiple characteristics is at present limited to various very special classes of examples [4],[9],[10],[8],[2].

As a test problem one might study operators

$$L = \sum_j X_j^2 + i \sum_{j,k} a_{j,k}[X_j, X_k]$$

where $\{X_j\}$ is a finite collection of $C^\infty$ vector fields satisfying the bracket hypothesis of Hörmander, and each coefficient $a_{j,k}$ is $C^\infty$ and real valued. Examples previously analyzed suggest that for fixed $\{X_j\}$, $L$ should be locally solvable for generic coefficients $a_{j,k}$. Most work to date has concerned examples of operators possessing high degrees of symmetry and depending on finitely many parameters, but recently the following examples have been analyzed [2]: $L = X^2 + Y^2 + ia(x)[X,Y]$ in $\mathbf{R}^3 = \{(x,y,t)\}$, where $X = \partial_x$ and $Y = \partial_y + x\partial_t$, and the coefficient $a$ depends on the coordinate $x$ alone but is otherwise arbitrary. It was shown that for real analytic $a$, $L$ is locally solvable at 0 unless $a$ is a constant function taking one of the values $\{\pm 1, \pm 3, \pm 5 \dots\}$.

*Date*: December 11, 1995.

Christ's research is supported in part by NSF grant DMS–9306833 and at MSRI by NSF grant DMS–9022140. Karadzhov's research is supported in part by a Fulbright Foundation grant.

This suggests that nonsolvability is far rarer, in relative terms, in the class of operators with general variable coefficients than it is within the special classes of highly symmetric operators discussed by earlier authors.

The purpose of the present paper is to exhibit a class of examples in which the nonsolvable operators form not a discrete set, but instead an infinite-dimensional manifold.[1] Thus nonsolvability is quite not so rare an occurrence as the results of [2] suggest it to be.

In $\mathbf{R}^2$ with coordinates $(x, t)$, define $X = \partial_x$ and $Y = x^{m-1}\partial_t$, where it is always assumed that $m \geq 2$, and consider

$$L = X^2 + Y^2 + ia(x)[X, Y] \tag{1.1}$$

where $a \in C^\infty$ is real valued. For each real constant $b$ define $\mathcal{L}_b = X^2 + Y^2 + ib[X, Y]$, and define $\Sigma_m$ to be the set of all $b \in \mathbf{R}$ such that $\mathcal{L}_b$ is not locally solvable at $0$. Denote by $a^{(n)}$ the $n$-th derivative of $a$. Our main result is then as follows.

**Theorem 1.1.** *Suppose that $m \geq 2$, $a \in C^\infty$ is real valued and that $L$ is defined as in (1.1).*
- *Solvability of $L$ at $0$ depends only on the Taylor coefficients of $a$ at $0$.*
- *Each set $\Sigma_m$ is nonempty and discrete.*
- *If $a(0) \notin \Sigma_m$ then $L$ is locally solvable at $0$.*
- *For each even $m \geq 2$ and each $b \in \Sigma_m$ there exist polynomials $P_j \colon \mathbf{R} \mapsto \mathbf{R}$ such that for all coefficients satisfying $a(0) = b$, $L$ is **not** locally solvable at $0$ if and only if for every $j \geq 1$,*

$$a^{(2j)}(0) = P_j\big(a^{(1)}(0), \ldots a^{(2j-1)}(0)\big). \tag{1.2}$$

The principal and only novel conclusion is the fourth. Consider the set of all equivalence classes of germs of smooth functions $a$ at $0$, with $a$ equivalent to $b$ if they have the same Taylor expansions at $0$. Local solvability of an operator $L$ of the form (1.1) depends only on the equivalence class of the coefficient $a$. The fourth conclusion is that for even $m$ there is an infinite dimensional family of equivalence classes corresponding to nonsolvable operators.

Curiously, the situation is quite different for odd $m$.

**Theorem 1.2.** *Suppose that $m \geq 3$ is odd, $a \in C^\infty$ is real valued and that $L$ is defined as in (1.1).*
- *If $m = 3$ then $L$ is locally solvable at $0$ unless $a(0) \in \Sigma_3$ and $a^{(j)}(0) = 0$ for every $j \geq 1$.*

---

[1] Our examples are nonetheless locally solvable for generic coefficients.

- *If $m \geq 5$ is odd then for each $b \in \Sigma_m$ there exist polynomials $P_j \colon \mathbf{R} \mapsto \mathbf{R}$ and a set $E \subset \mathbf{N}$ of cardinality at most $m - 3$ such that for all coefficients satisfying $a(0) = b$, $L$ is **not** locally solvable at $0$ if and only if for every $1 \leq j \notin E$,*

$$a^{(j)}(0) = P_j\big(a^{(1)}(0), \ldots a^{(j-1)}(0)\big). \qquad (1.3)$$

There is a simple explanation for why there might be many more nonsolvable operators here than in the superficially similar three-dimensional situation studied in [2]. If $a(0) \in \Sigma_m$, then local solvability of $L$ is closely related to the asymptotic behavior, as $\varepsilon \to 0$, of the unique small eigenvalue $\lambda(\varepsilon)$ of the ordinary differential operator $B_\varepsilon = -\partial_y^2 + y^{2(m-1)} + (m-1)y^{m-2}a(\varepsilon y)$, as will be seen in Section 2; $\lambda(\varepsilon)$ has an asymptotic expansion in powers of $\varepsilon$, and nonsolvability is equivalent to the vanishing of all coefficients in this expansion. Thus we have a type of inverse spectral problem. The situation in [2] is analogous, except that $\lambda$ depends there on two parameters while $a$ still depends on only one, so the inverse problem is formally overdetermined and hence it is plausible that $\lambda$ vanishes to infinite order at $00$ only for very few coefficients $a$.

These considerations fail to distinguish, however, between odd and even parameters $m$. When $m$ is even there is an extra symmetry in the problem that is absent for odd $m$: $B_\varepsilon$ is unitarily equivalent to $B_{-\varepsilon}$, hence $\lambda(\varepsilon)$ is an even function. Thus the inverse problem is highly underdetermined; the situation is quite analogous to the classical inverse spectral problem for $-\partial_x^2 + q(x)$ on $[-1,1]$ [13]. For odd $m$ the symmetry is broken, and the inverse problem has a unique solution for $m = 3$ and is underdetermined by at most finitely many parameters for larger odd $m$.

One corollary of our theorem is that the nonsolvable operators are not all reducible, via obvious symmetries, to examples with constant or very simple coefficients.

## 2. Solvability and eigenvalue asymptotics

It will always be assumed that $m \geq 2$. The notation $\|\cdot\|$, with no subscript, denotes always the norm in $L^2(\mathbf{R})$ with respect to Lebesgue measure.

Fix $\eta \in C_0^\infty(\mathbf{R})$ supported in an open set where $a$ is defined, satisfying $\eta \equiv 1$ in some smaller neighborhood of $0$. Redefine $b(x) = [a(x) - a(0)]\eta(x)$ and consider

$$L = \partial_x^2 + x^{2(m-1)}\partial_t^2 + i(m-1)x^{m-2}(a(0) + b(x))\partial_t. \qquad (2.1)$$

This agrees with our original operator $L$ in a neighborhood of $0$, so it suffices to study its local solvability instead.

Applying the partial Fourier transform with respect to the coordinate $t$ leads to the one parameter family of ordinary differential operators

$$A_\tau = \partial_x^2 - x^{2(m-1)}\tau^2 - (m-1)\tau x^{m-2}(a(0) + b(x)), \qquad \tau \in \mathbf{R}.$$

**Lemma 2.1.** *Suppose there exist $c, M \in \mathbf{R}^+$ and a neighborhood $V$ of $0$ such that for all sufficiently large $|\tau|$,*

$$\|f\|_{L^2} \leq C|\tau|^M \|Lf\|_{L^2} \tag{2.2}$$

*for all $f \in C_0^2(V)$. Then $L$ is locally solvable at $0$.*

This is proved in a more general context in the second paragraph of Section 2 of [2] with $L$ replaced by its transpose on the right hand side of (2.2). Replacing the transpose by the adjoint results in an equivalent inequality, by Taking complex conjugates, and in the present case $L$ equals its adjoint. □

Assume that $\tau \neq 0$. Substituting $x = \varepsilon y$, where $\varepsilon = |\tau|^{-1/m}$, $A_\tau$ is unitarily equivalent either to $-\tau^{2/m} B_\varepsilon^+$ or to $-\tau^{2/m} B_\varepsilon^-$, depending on the sign of $\tau$, where

$$B_\varepsilon^\pm = -\partial_y^2 + y^{2(m-1)} \pm (m-1)y^{m-2}(a(0) + b(\varepsilon y)). \tag{2.3}$$

Each $B_\varepsilon^\pm$ is an essentially selfadjoint[2] operator on $L^2(\mathbf{R})$, with compact resolvent, whose spectrum consists of a discrete sequence of eigenvalues, all of multiplicity one, tending to $+\infty$. Denote those eigenvalues by $\mu_0^\pm(\varepsilon) < \mu_1^\pm(\varepsilon) < \ldots$. Then by Lemma 2.1, in order to prove local solvability of $L$ it suffices to prove that for both choices of the $\pm$ sign,

$$\min_j |\mu_j^\pm(\varepsilon)| \geq c\varepsilon^M \qquad \text{as } \varepsilon \to 0^+, \tag{2.4}$$

for some $c, M \in \mathbf{R}^+$.

Define $\Sigma_m^+$ to be the set of all $\alpha \in \mathbf{R}$ such that $0$ is an eigenvalue of

$$\mathcal{L}_\alpha = -\partial_y^2 + y^{2(m-1)} + (m-1)\alpha y^{m-2}.$$

Define $\Sigma_m = \Sigma_m^+ \cup (-\Sigma_m^+) = \{\alpha : \alpha \in \Sigma_m^+ \text{ or } -\alpha \in \Sigma_m^+.\}$.

**Lemma 2.2.** *$\Sigma_m^+$ is a nonempty, discrete subset of $\mathbf{R}$.*

*Proof.* From the identity

$$\langle \mathcal{L}_0 f, f \rangle = \int_\mathbf{R} \left( |\partial_y f|^2 + |y^{m-1} f|^2 \right) dy \tag{2.5}$$

it follows that the domain of $\mathcal{L}_\alpha$ is independent of $\alpha$, and hence that $\alpha \mapsto \mathcal{L}_\alpha$ is an analytic family of unbounded operators in the sense of Kato [6]. Since the eigenvalues of $\mathcal{L}_\alpha$ are all simple, each eigenvalue therefore extends to an entire holomorphic function of $\alpha \in \mathbf{C}$. Thus either $\Sigma_m^+$ is discrete, or 0 is an eigenvalue of $\mathcal{L}_\alpha$ for every $\alpha \in \mathbf{C}$.

---

[2]See the proof of Lemma 3 of [3].

(2.5) implies that $\|\mathcal{L}_0 f\| \geq c\|f\|$ for some $c > 0$, for all $f$ in the domain of $\mathcal{L}_0$, since the right hand side of (2.5) is bounded below by some positive constant times $\|f\|^2$ by an elementary argument. Thus 0 is not an eigenvalue of $\mathcal{L}_0$, so $\Sigma_m^+$ is discrete.

$\mathcal{L}_\alpha$ is essentially selfadjoint for each $\alpha \in \mathbf{R}$, and $\mathcal{L}_0$ is strictly positive. In order to show that 0 is an eigenvalue of some $\mathcal{L}_\alpha$, it suffices by the intermediate value theorem and continuous dependence of the eigenvalues on $\alpha$ to show that there exist $f \in C_0^2(\mathbf{R})$ and $\beta \in \mathbf{R}$ such that $\langle \mathcal{L}_\beta f, f \rangle < 0$. But fixing any $f$ that is supported in $\mathbf{R}^+$ and does not vanish identically, this inequality holds for all sufficiently negative $\beta$. □

Most of the information we require on the operators $B_\varepsilon^\pm$ is either implicit in [2], or follows directly from small modifications of that analysis, often with simplifications. Therefore we will not repeat the details of the proofs,[3] but will state the results needed and refer the reader to the relevant passages in [2].

**Lemma 2.3.** *If $a(0) \notin \Sigma_m$ then*

$$\min_j |\mu_j^\pm(\varepsilon)| \geq c > 0 \qquad as\ \varepsilon \to 0. \tag{2.6}$$

*Therefore if $a(0) \notin \Sigma_m$ then $L$ is locally solvable at 0.*

The proof of (2.6) is a simplification of the proofs of Lemma 3.2 and 3.3 of [2]. Actually $L$ satisfies maximal and subelliptic estimates in this case, as follows from general theory [12], or from elementary arguments as in [3]. □

We now restrict attention to $B_\varepsilon^+$; all reasoning will apply equally to $B_\varepsilon^-$ after minor changes of notation. To simplify the notation we then drop the superscript $\pm$, writing for all $\varepsilon \in \mathbf{R}$

$$B_\varepsilon = -\partial_y^2 + y^{2(m-1)} + (m-1)y^{m-2}(a(0) + b(\varepsilon y)).$$

**Lemma 2.4.** *If $a(0) \in \Sigma_m$ then there exists $\theta > 0$ such that for all sufficiently small $\varepsilon \in \mathbf{R}$, $B_\varepsilon$ has precisely one eigenvalue in $[-\theta, \theta]$. This eigenvalue tends to 0 as $\varepsilon \to 0$.*

---

[3]Alternatively, Lemmas 2.3 and 2.4 may be proved quite simply by observing that for all small $\varepsilon > 0$, $(1-\varepsilon)B_0 - C\varepsilon \leq B_\varepsilon \leq (1+\varepsilon)B_0 + C\varepsilon$ for some fixed positive constant $C$. This simple argument does not work in all cases in [2]. Lemmas 2.5 and 2.7 may be avoided by appealing to the theory of globally elliptic operators.

For a proof see Lemma 3.4 of [2]. □

Assume until further notice that $a(0) \in \Sigma_m$. Denote by $\lambda(\varepsilon)$ the unique small eigenvalue of $B_\varepsilon$, and fix a real valued eigenfunction $\psi_0$ of $B_0$ satisfying $\|\psi_0\| = 1$, associated to the eigenvalue $\lambda(0) = 0$.

**Lemma 2.5.** *There exists $\delta > 0$ such that*
$$\int_{\mathbf{R}} \left(|\partial_y \psi_0|^2 + \psi_0^2\right) e^{\delta |y|} \, dy < \infty. \tag{2.7}$$
*Moreover there exists $C < \infty$ such that for all $0 \leq r \leq \delta$ and all $f \in L^2$ orthogonal to $\psi_0$,*
$$\int_{\mathbf{R}} (B_0^{-1} f)^2 \, e^{r|y|} \, dy \leq C \int_{\mathbf{R}} f^2 \, e^{r|y|} \, dy. \tag{2.8}$$

See the proofs of Lemmas 3.2 and 3.5 of [2]. □

By $B_0^{-1}$ we mean the following operator. $B_0$ has a one dimensional kernel and cokernel, both spanned by $\psi_0$. $B_0^{-1}$ is the unique bounded linear operator on $L^2(\mathbf{R})$ satisfying $B_0^{-1} g = 0$ for any multiple $g$ of $\psi_0$, $B_0 B_0^{-1} g = g$ for every $g \perp \psi_0$, and $B_0^{-1} g \perp \psi_0$ for every $g \in L^2$.

**Lemma 2.6.** *There exist coefficients $\Lambda_j \in \mathbf{R}$ such that for every $N$,*
$$\lambda(\varepsilon) = \sum_{j=0}^{N} \Lambda_j \varepsilon^j + O(\varepsilon^{N+1}) \qquad \text{as } \varepsilon \to 0. \tag{2.9}$$

See the proof of Lemma 3.7 of [2]. □

For our purpose it is necessary to recall from [2] the algebraic formalism determining the coefficients $\Lambda_j$. Expand
$$B_\varepsilon \sim \sum_{j=0}^{\infty} \beta_j \varepsilon^j$$
as a formal operator valued Taylor series about $\varepsilon = 0$, where
$$\beta_0 = -\partial_y^2 + y^{2(m-1)} + (m-1) y^{m-2} a(0) = B_0$$
and
$$\beta_j = y^{m-2+j} \cdot (m-1) a^{(j)}(0)/j! \qquad \text{for all } j \geq 1.$$
Writing an eigenfunction with eigenvalue $\lambda(\varepsilon)$ formally as $\varphi(\varepsilon) = \sum_j \varphi_j \varepsilon^j$ with $\varphi_0 = \psi_0$, then equating coefficients of like powers of $\varepsilon$ in the formal power series equation $B_\varepsilon \varphi(\varepsilon) \sim \lambda(\varepsilon) \varphi(\varepsilon)$, gives the following relation.
$$B_0 \varphi_n = -\sum_{j=1}^{n} (\beta_j - \Lambda_j) \varphi_{n-j}. \tag{2.10}$$

In order that (2.10) should have a solution in $L^2$ it is necessary and sufficient that the right hand side be orthogonal to $\psi_0$, which spans the cokernel of $B_0$. Since $\|\psi_0\| = 1$, this orthogonality condition is equivalent to

$$\Lambda_n = \langle \beta_n \psi_0, \psi_0 \rangle + \sum_{j=1}^{n-1} \langle (\beta_j - \Lambda_j) \varphi_{n-j}, \psi_0 \rangle \tag{2.11}$$

where the summation on the right is empty for $n = 1$. If (2.11) is satisfied and the right hand side of (2.10) belongs to $L^2$ there exists a unique solution $\varphi_n$ orthogonal to the nullspace of $B_0$, and we take always this distinguished solution, so that $\varphi_n$ is uniquely determined. (2.11) may thus be restated as

$$\Lambda_n = \langle \beta_n \psi_0, \psi_0 \rangle + \sum_{j=1}^{n-1} \langle \beta_j \varphi_{n-j}, \psi_0 \rangle \tag{2.12}$$

since the terms omitted are zero.

Once $\varphi_0, \ldots \varphi_n$ and $\Lambda_0, \ldots \Lambda_n$ have been determined, $\Lambda_{n+1}$ is defined by (2.11). Then (2.10) is solvable and determines $\varphi_{n+1}$. Thus all $\varphi_j$ and $\Lambda_j$ are determined by induction on $j$. Lemma 2.5 and an induction argument demonstrate that for each $n$, $\varphi_n(y) = O(\exp(-r_n|y|))$ for some $r_n > 0$, so that the right hand side of (2.10) does belong to $L^2(\mathbf{R})$, and the right hand side of (2.11) is well defined.

**Lemma 2.7.** *Each function $\varphi_n$ belongs to the Schwartz class $\mathcal{S}$.*

This is proved by induction on $n$. Thus it is given that the right hand side of (2.10) belongs to $\mathcal{S}$, and from Lemma 2.5 we know that $\varphi_n$ and $\partial_y \varphi_n$ decay exponentially, in the $L^2$ sense. Since $B_0$ has polynomial coefficients, a straightforward bootstrapping argument shows that if $f, \partial_y f$, and $B_0 f$ decay exponentially, in the $L^2$ sense, then so does $\partial_y^k f$, for all $k$. (2.10) and induction on $n$ then imply that each derivative $\partial_y^k \varphi_n$ decays exponentially, in the $L^2$ sense. $\square$

**Lemma 2.8.** *If $m$ is even then for all coefficients $a$, $\Lambda_n = 0$ for every odd $n \geq 1$.*

*Proof.* For even $m$, $B_\varepsilon$ is unitarily equivalent to $B_{-\varepsilon}$ for all $\varepsilon$, as is seen by the change of variables $y \mapsto -y$. Thus $\lambda(\varepsilon) \equiv \lambda(-\varepsilon)$, and hence $\Lambda_n$ must vanish for all odd $n$. $\square$

**Lemma 2.9.** *If $a(0) \in \Sigma_m$ remains fixed then each $\Lambda_n$ is a polynomial in the quantities $a^{(1)}(0), \ldots a^{(n)}(0)$.*

The coefficients of this polynomial depend only on $m$ and on $a(0)$.

*Proof.* In order to prove this by induction on $n$, we claim further that for each $n$ there exist finitely many functions $\varphi_{n,i} \in L^2(\mathbf{R}, \exp(r_n|y|)dy)$ and polynomials $q_{n,i}[a^{(1)}(0) \ldots a^{(n)}(0)]$ such that

$$\varphi_n = \sum_i q_{n,i}[a^{(1)}(0), \ldots a^{(n)}(0)]\varphi_{n,i} \tag{2.13}$$

for all $a$ with prescribed value $a(0)$. Supposing these claims to hold for all $n \leq N$, it follows from (2.11) that $\Lambda_{N+1}$ has the desired form.[4] By (2.10),

$$\varphi_{N+1} = -\sum_{j=1}^{N+1} \sum_i q_{N+1-j,i} B_0^{-1}(\beta_j - \Lambda_j)\varphi_{N+1-j,i}$$

then also takes the desired form. □

From (2.10) it is now immediate that for any $m, n$, if $a(0) \in \Sigma_m$ remains fixed, then

$$\Lambda_n = c_n \cdot a^{(n)}(0) - Q_n(a^{(1)}(0) \ldots a^{(n-1)}(0)) \tag{2.14}$$

for certain polynomials $Q_n$ which depend also on $m$ and on $a(0)$, with

$$c_n = [(m-1)/n!]\int_{\mathbf{R}} \psi_0^2(y) y^{m-2+n}\, dy. \tag{2.15}$$

In particular

$$c_n \neq 0 \text{ for all pairs } (m, n) \text{ such that } m+n \text{ is even.}$$

In the special case $m = 2$, $\Sigma_2^+ = \{-1, -3, -5, \ldots\}$ and all the $Q_n, c_n$ can in principle be computed using Hermite polynomials, since the product of $y$ with any Hermite polynomial can be expressed explicitly as a linear combination of two Hermite polynomials.

The same analysis applies to the operators $B_\varepsilon^-$, with $\Sigma_m^+$ replaced by $-\Sigma_m^+$. For even $m$, because the function $y^{m-2}$ is nonnegative, $\Sigma_m^+$ is clearly contained in $(-\infty, 0)$, hence $-\Sigma_m^+ \cap \Sigma_m^+ = \emptyset$. When $a(0) \in -\Sigma_m^+$ we define $c_n, Q_n$ in terms of $B_\varepsilon^-$ by the same procedure as above. When $m$ is odd, $\Sigma_m^+ = -\Sigma_m^+ = \Sigma_m$, since the change of variables $y \mapsto -y$ leaves $-\partial_y^2 + y^{2(m-2)}$ invariant but reverses the sign of $y^{m-2}$. Thus $B_\varepsilon^+$ is unitarily equivalent to $B_{-\varepsilon}^-$, so the corresponding small eigenvalues $\lambda^\pm(\varepsilon)$ satisfy $\lambda^+(\varepsilon) = \lambda^-(-\varepsilon)$ for all $\varepsilon$, and hence the coefficients $\Lambda_j^\pm$ in their asymptotic expansions are related by $\Lambda_j^+ = (-1)^j \Lambda_j^-$. We have thus proved the following result.

---

[4] Nonlinear dependence stems from the occurrence of $\beta_j \varphi_{n-j}$.

**Proposition 2.10.** *Let $m \geq 2$ and suppose that $a(0) \in \Sigma_m$. Suppose that there exists $j \geq 1$ such that $m - j$ is even and*

$$a^{(j)}(0) \neq Q_j(a^{(1)}(0), \ldots a^{(j-1)}(0)).$$

*Then*

$$\min_i |\mu_i^{\pm}(\varepsilon)| \geq c|\varepsilon|^j \qquad \text{as } \varepsilon \to 0$$

*for both choices of the $\pm$ sign. In particular, $L$ is locally solvable at $0$.*

The case of odd $m$ will be taken up in Section 4.

## 3. Nonsolvability

The purpose of this section is to establish Proposition 3.1, a converse to Lemma 2.10. For even $m$ we denote by $\lambda(\varepsilon)$ the small eigenvalue of $B_\varepsilon^+$ if $a(0) \in \Sigma_m^+$, and of $B_\varepsilon^-$ if $a(0) \in -\Sigma_m^+$.

**Proposition 3.1.** *Suppose that $a(0) \in \Sigma_m$ and that $\lambda(\varepsilon) = O(\varepsilon^N)$ as $\varepsilon \to 0$, for all $N$. Then $L$ is not locally solvable at $0$.*

We may suppose without loss of generality that $a(0) \in \Sigma_m^+$, by replacing $a$ by $-a$ if necessary. We are thus given that $\Lambda_j = 0$ for all $j \geq 0$, and that formally $B_\varepsilon(\sum \varepsilon^j \varphi_j) = O(\varepsilon^N)$ for all $N$.

Set

$$\Phi_\varepsilon = \Phi_\varepsilon^{(A)} = \sum_{j=0}^{A} \varepsilon^j \varphi_j,$$

where $A$ is a large positive integer to be chosen later. Formally $B_\varepsilon \Phi_\varepsilon = O(\varepsilon^{A+1})$. Since each $\varphi_j$ is a Schwartz function, and all $\Lambda_j$ vanish, we have more precisely

$$|\partial_y^k B_\varepsilon \Phi_\varepsilon(y)| \leq C_{k,M} |\varepsilon|^A (1 + |y|)^{-M} \tag{3.1}$$

for all $k, M$, as follows from (2.13) and the fact that the coefficient of $\varepsilon^j$, in the formal expansion for $B_\varepsilon \Phi_\varepsilon$ about $\varepsilon = 0$, vanishes for all $j \leq A$. Define

$$F_\tau(x) = F_\tau^{(A)}(x) = \Phi_{\tau^{-1/2}}(\tau^{1/2} x)$$

for $\tau \in \mathbf{R}^+$. Given any $N < \infty$, (3.1) asserts that

$$|\partial_x^\alpha A_\tau F_\tau(x)| \leq C_N \tau^{-N} (1 + |\tau^{1/2} x|)^{-N} \tag{3.2}$$

for all $x \in \mathbf{R}$, all $\tau \geq 1$ and all $0 \leq \alpha \leq N$, provided that $A$ is chosen to be sufficiently large.

Fix $\eta \in C_0^\infty(\mathbf{R})$ supported in $(4^{-1}, 4)$ and $\equiv 1$ on $[2^{-1}, 2]$. For large $\lambda \in \mathbf{R}^+$ define

$$G_\lambda(x, t) = \int_\mathbf{R} e^{it\tau} \eta(\tau/\lambda) F_\tau(x) \, d\tau$$

$$= \int e^{it\tau} \eta(\tau/\lambda) \sum_{j=0}^A \tau^{-j/2} \varphi_j(\tau^{1/2} x) \, d\tau.$$

Since each $\varphi_j$ is a Schwartz function,

$$|\partial_x^\alpha F_\tau(x)| \leq C \lambda^{\alpha/2} (1 + \lambda^{1/2} |x|)^{-M}$$

for any $\alpha, M$, for all $\tau \in (4^{-1}\lambda, 4\lambda)$, uniformly in $\lambda$. Writing

$$\partial_t^\beta \partial_x^\alpha G_\lambda(x, t) = \int e^{it\tau} (i\tau)^\beta \eta(\tau/\lambda) \partial_x^\alpha F_\tau(x) \, d\tau, \tag{3.3}$$

it follows that

$$|\partial_t^\beta \partial_x^\alpha G_\lambda(x, t)| \leq C \lambda^{1+\beta+\alpha/2} (1 + \lambda^{1/2} |x|)^{-M} \tag{3.4}$$

for all $\alpha, \beta, M$, for some $C_{\alpha,\beta,M}$. Integrating by parts $n$ times with respect to $\tau$ in (3.3) gives

$$\partial_t^\beta \partial_x^\alpha G_\lambda(x, t) = \pm i^{n+\beta} t^{-n} \int e^{it\tau} \partial_\tau^n \left[ \tau^\beta \eta(\tau/\lambda) F_\tau(x) \right] d\tau.$$

Thus

$$|\partial_t^\beta \partial_x^\alpha G_\lambda(x, t)| \leq C \lambda^{1+\beta+\alpha/2} (1 + \lambda^{1/2} |x|)^{-M} (\lambda |t|)^{-n} \tag{3.5}$$

for any $\alpha, \beta, M, m$. This bound will be used in the weaker form

$$|\partial_t^\beta \partial_x^\alpha G_\lambda(x, t)| \leq C_M \lambda^{-M} \text{ whenever } |x| \geq \lambda^{-1/4} \text{ or } |t| \geq \lambda^{-1/2}, \tag{3.6}$$

as follows from (3.5) by choosing $M, n$ to be sufficiently large.

Applying Plancherel's theorem in the coordinate $t$ and recalling that $\varphi_0 = \psi_0 \neq 0$, we obtain a lower bound

$$\|G_\lambda\|_{L^2(\mathbf{R})}^2 \geq c \lambda^{1/2}$$

for all sufficiently large $\lambda$, for some $c > 0$. So for large $\lambda$ there certainly exists $z = z(\lambda) \in \mathbf{R}^2$ such that $|z_1| \leq \lambda^{-1/4}$ and $|z_2| \leq \lambda^{-1/2}$, satisfying $|G_\lambda(z)| \geq 1$.

By (3.5), the $L^2$ norm of $G_\lambda$ on the region where $|x| \geq \lambda^{-1/4}$ or $|t| \geq \lambda^{-1/2}$ is $O(\lambda^{-M})$ for all $M$. By (3.4), $\nabla G_\lambda = O(\lambda^2)$ in the supremum norm. Therefore $|G_\lambda(z') - G_\lambda(z)| \leq 1/2$ for all $|z' - z| \leq \lambda^{-3}$, for all sufficiently large $\lambda$. Fix $h \in C_0^\infty(\mathbf{R}^2)$ satisfying $\int h = 1$, and set

$$h_\lambda(w) = h(\lambda^4 (w - z(\lambda)))$$

for $w \in \mathbf{R}^2$. Then $h_\lambda$ is supported where $|w - z(\lambda)| \leq C\lambda^{-4} \ll \lambda^{-3}$, so

$$\left|\int_{\mathbf{R}^2} G_\lambda \overline{h_\lambda}\right| \geq c\lambda^{-8} \tag{3.7}$$

for some $c > 0$, for all sufficiently large $\lambda$. On the other hand,

$$\|h_\lambda\|_{C^M} \leq C_M \lambda^{4M}$$

for any $M$.

Fix yet another cutoff function $\zeta \in C_0^\infty(\mathbf{R}^2)$ satisfying $\zeta(x,t) \equiv 1$ for $|(x,t)| \leq 2$ and define

$$g_\lambda(x,t) = G_\lambda(x,t)\zeta(\lambda^{1/4}x, \lambda^{1/2}t).$$

From (3.2) together with (3.6) it follows that

$$Lg_\lambda(x,t) = O(\lambda^{-N})$$

in any $C^B$ norm as $\lambda \to \infty$, where $N$ may be made arbitrarily large by choosing $A$ to be sufficiently large. Moreover $g_\lambda$ is supported in the region where $|(x,t)| \leq C\lambda^{-1/4}$, which shrinks to 0 as $\lambda \to \infty$.

If $L$ were locally solvable at 0 there would exist [5] $B < \infty$ such that

$$\left|\int g\overline{h}\right| \leq B\|h\|_{C^B}\|L^*g\|_{C^B} \tag{3.8}$$

for any $h, g \in C_0^\infty$ supported in $\{|(x,t)| \leq B^{-1}\}$. The preceding analysis shows that no such inequality is valid. Fix any constant $B$. Taking $g = g_\lambda, h = h_\lambda$ and noting that $L = L^*$, $\|L^*g\|_{C^B} = O(\lambda^{-N})$ where $N$ may be taken to be arbitrarily large by choosing $A$ to be sufficiently large. $h_\lambda$ is not also small, but $\|h_\lambda\|_{C^B} \leq C\lambda^{4B}$, so the right hand side in (3.8) is $O(\lambda^{-N+4B}) = O(\lambda^{-N/2})$ for any large preassigned $N$. On the other hand, $|\int g_\lambda \overline{h_\lambda}| \geq c\lambda^{-8}$. Thus given $B$, $N$ may be chosen so that (3.8) becomes false for all sufficiently large $\lambda$. Thus $L$ is not locally solvable at 0. $\square$

## 4. The case of odd $m$

Suppose that $m \geq 3$ is odd, and that $a(0) \in \Sigma_m$. Then $\Sigma_m^+ = \Sigma_m^-$, and we begin by discussing $B_\varepsilon^+$ and the coefficients $\Lambda_j = \Lambda_j^+$.

With the notation of Lemma 2.6 and of the discussion following it, consider the moments

$$\mu_j = \int_{\mathbf{R}} x^j \psi^2(x)\, dx$$

for each integer $j \geq 0$, where $\psi = \psi_0$. These depend of course on $m$ and on $a(0)$. Define $\mu_j = 0$ for $j < 0$.

**Lemma 4.1.** *For each $j \geq -2$,*

$$(2m + 2j + 2)\mu_{2m+j-1} + (m + 2j + 2)\alpha\mu_{m+j-1} = \tfrac{1}{2}j(j+1)(j+2)\mu_{j-1}. \tag{4.1}$$

*Proof.* We have $B_0\psi = 0$, that is,

$$\psi'' = (x^{2m-2} + \alpha x^{m-2})\psi. \tag{4.2}$$

Multiplying both sides of this equation by $x^{j+1}\psi$ and integrating by parts gives, for all $j \geq -1$,

$$\int (x^{2m+j-1} + \alpha x^{m+j-1})\psi^2 = -\int x^{j+1}(\psi')^2 + \tfrac{j(j+1)}{2}\int x^{j-1}\psi^2, \tag{4.3}$$

where all integrals are taken over $\mathbf{R}$ with respect to $dx$. The final integral is not well defined when $j \leq 0$, but the identity remains valid if this integral is interpreted to be $\mu_{j-1} = 0$. Similarly, for any $j \geq -2$, multiplying both sides of (4.2) by $x^{j+2}\psi'$ and integrating by parts in all terms gives

$$\int \big((2m+j)x^{2m+j-1} + (m+j)\alpha x^{m+j-1}\big)\psi^2 = (j+2)\int x^{j+1}(\psi')^2, \tag{4.4}$$

interpreting the right hand side to be zero when $j = -2$. Taking an appropriate linear combination of (4.3) and (4.4) yields (4.1) for each $j \geq -1$. When $j = -2$, (4.4) alone suffices to give (4.1). $\square$

Note that $\mu_j > 0$ for all even $j$.

**Lemma 4.2.** *If $m = 3$ then for every $j \geq 0$, $\mu_j \neq 0$. If $m > 3$ is odd then there are at most $m - 3$ indices $j$ satisfying $\mu_j = 0$.*

*Proof.* The recurrence relation (4.1) simplifies to

$$(2m + 2j + 2)\mu_{2m+j-1} + (m + 2j + 2)\alpha\mu_{m+j-1} = 0 \qquad \text{for } j = -2, -1, 0. \tag{4.5}$$

$0 \notin \Sigma_m$ since $-\partial_y^2 + y^{2(m-1)}$ is a strictly positive operator, so $\alpha \neq 0$. Since $m$ is odd and $\mu_j \neq 0$ for even $j \geq 0$, (4.5) then implies that $\mu_i \neq 0$ for all $i \in \{m-3, m-2, m-1, 2m-3, 2m-2, 2m-1\}$. Moreover when $m = 3$, $\mu_{i+m}$ and $\mu_i$ have opposite signs for $i = 0, 1, 2$ if $\alpha > 0$, and have the same signs if $\alpha < 0$.

First suppose that $m = 3$ and $\alpha < 0$. Then $\mu_i > 0$ for all $0 \leq i < 6$, and the recurrence relation (4.1) then implies that $\mu_j > 0$ for all $j \geq 0$. If $m = 3$ but $\alpha > 0$, then $\mu_i > 0$ for $i = 0, 2, 4$ and $\mu_i < 0$ for $i = 1, 3, 5$, and the recurrence relation then implies by induction on $j$ that $\mu_j > 0$ for all even $j \geq 0$ and $\mu_j < 0$ for all odd $j \geq 1$.

For odd $m > 3$ this reasoning breaks down to a degree, because (4.5) still gives only 3 initial relations, while (4.1) decomposes into $m$ independent three term recurrence relations. The same reasoning as for $m = 3$ applies to 3 of these, but $m - 3$ are left.

Fix some $r \in \{0, 1, 2, \ldots m-1\}$ and consider the sequence of coefficients $b_i = \mu_{r+im}$. Then $b_i > 0$ whenever $r + im$ is even, that is, for every second index $i$. The recurrence relation (4.1) is a three term recurrence relation for the sequence $b_i$, for each fixed $r$. We claim that for fixed $r$, there can be at most one index $k$ for which $b_k = 0$. Indeed, if $k$ is such an index, then setting $j = km + r + 1$ in (4.1) yields

$$(2m + 2(km + r + 1) + 2)b_{k+2} + (m + 2(km + r + 1) + 2)\alpha b_{k+1} = 0.$$

This analogue of (4.5) may then be used as an initial condition and combined with the recurrence relation (4.1) just as in the case $m = 3$, to deduce that $b_i \neq 0$ for all nonnegative $i \neq k$. □

To prove Theorem 1.2 for $m = 3$, suppose that $L$ is not locally solvable at 0. Then $a(0) \in \Sigma_m$ and all coefficients $\Lambda_j^+ = (-1)^j \Lambda_j^-$ vanish. Write $\Lambda_j = \Lambda_j^+$. We have

$$\Lambda_1 = \langle \beta_1 \psi_0, \psi_0 \rangle = a^{(1)}(0)\mu_1.$$

Since $\mu_1 \neq 0$, $a^{(1)}(0)$ must vanish. Suppose, by induction on $n$, that $a^{(j)}(0)$ has been proved to vanish for all $1 \leq j < n$. Thus $\beta_j = 0$ for all $1 \leq j < n$, and hence (2.10) gives $B_0 \varphi_n = 0$, so $\varphi_n = 0$ since it is by definition orthogonal to the nullspace of $B_0$. (2.11) then simplifies to

$$\Lambda_n = \langle \beta_n \psi_0, \psi_0 \rangle = (m-1)a^{(n)}(0)\mu_n/n! .$$

Since $\Lambda_n$ is assumed to vanish and $\mu_n \neq 0$, this forces $a^{(n)}(0) = 0$.

The reasoning for the case of odd $m \geq 5$ is similar, except that if $\mu_n$ does vanish then no constraint can be deduced on $a^{(n)}(0)$.

MICHAEL CHRIST, UNIVERSITY OF CALIFORNIA, LOS ANGELES, LOS ANGELES, CA. 90095-1555
  *E-mail address*: christ@math.ucla.edu

G. E. KARADZHOV, INSTITUTE OF MATHEMATICS, BULGARIAN ACADEMY OF SCIENCE, P.O. BOX 373, 1090 SOFIA, BULGARIA
  *E-mail address*: difeq@bgearn.bitnet

DETLEF MÜLLER, MATHEMATISCHES SEMINAR, CHRISTIAN ALBRECHTS-UNIVERSITÄT KIEL, LUDEWIG-MEYN-STR.4, D-24098 KIEL, GERMANY
  *E-mail address*: mueller@math.uni-kiel.de